\documentclass[aop,MSNbibl,seceqn,citesort,dvips]{arximspdf}
\usepackage{mathbh}

\doi{10.1214/11-AOP691}
\volume{41}
\issue{3B}
\pubyear{2013}
\firstpage{1806}
\lastpage{1830}

\makeatletter

\newtheorem{lemma}{Lemma}[section]
\newtheorem{theorem}{Theorem}[section]
\newproclaim{remark}{Remark}[section]

\newcommand{\D}{\Delta}
\newcommand{\e}{\varepsilon}
\newcommand{\g}{\gamma}
\newcommand{\s}{\sigma}

\newcommand{\p}{\partial}

\newcommand{\bN}{{\mathbb N}}
\newcommand{\bR}{{\mathbb R}}
\newcommand{\bP}{{\mathbb P}}
\newcommand{\bE}{{\mathbb E}}
\newcommand{\bZ}{{\mathbb Z}}

\newcommand{\cA}{{\mathcal A}}
\newcommand{\cD}{{\mathcal D}}
\newcommand{\cP}{{\mathcal P}}

\newcommand{\bL}{{\mathbf L}}
\newcommand{\bff}{{\mathbf f}}

\makeatother

\begin{document}
\begin{frontmatter}

\title{Gelation for Marcus--Lushnikov process}
\runtitle{Gelation for Marcus--Lushnikov process}

\begin{aug}
\author[A]{\fnms{Fraydoun} \snm{Rezakhanlou}\corref{}\thanksref{t1,t2}\ead[label=e1]{rezakhan@math.berkeley.edu}}
\runauthor{F. Rezakhanlou}
\affiliation{University of California, Berkeley}
\address[A]{Department of Mathematics\\
University of California, Berkeley\\
Berkeley, California 94720-3840\\
USA\\
\printead{e1}} 
\end{aug}

\thankstext{t1}{Supported in part by NSF Grant DMS-07-07890.}

\thankstext{t2}{Supported by IHES and Universite de Paris 12.}

\received{\smonth{1} \syear{2011}}
\revised{\smonth{7} \syear{2011}}

%
\begin{abstract}
The Marcus--Lushnikov process is a simple mean field model of
coagulating particles that converges to the homogeneous Smoluchowski
equation in the large mass limit. If the coagulation rates grow
sufficiently fast as the size of particles get large, giant particles
emerge in finite time. This is known as gelation, and such particles
are known as gels. Gelation comes in different flavors: simple,
instantaneous and complete. In the case of an instantaneous gelation,
giant particles are formed in a very short time. If all particles
coagulate to form a single particle in a time interval that stays
bounded as total mass gets large, then we have a complete gelation. In
this article, we describe conditions which guarantee any of the three
possible gelations with explicit bounds on the size of gels and the
time of their creations.
\end{abstract}

%

\end{frontmatter}

\section{Introduction}
\label{sec1}

The Smoluchowski equation is a coupled system of differential equations
that describes the evolving densities
(or concentrations) of a
system of particles (or clusters)
that are prone to coagulate in pairs. A~sequence of functions
$f_n\dvtx [0,\infty)\to[0,\infty)$, $n \in\bN$, is a solution of the (discrete
and homogeneous) Smoluchowski equation (SE) if it satisfies
%
\begin{equation}
\label{eq11}
\frac{d}{d t} f_n(t) = Q_n(f)(t)
\end{equation}
with $Q_n = Q_n^{+} - Q_n^{-}$, where
\begin{eqnarray*}
Q_n^{+}(f)(t) &=& \frac12\sum_{m=1}^{n-1}
\alpha(m,n-m)f_m(t)f_{n-m}(t),\\
Q_n^{-}(f)(t) &=& \sum_{m=1}^{\infty} \alpha(n,m)f_n(t)f_m(t).
\end{eqnarray*}

The function $f_n$ represents the density of particles of size $n$, and
the symmetric function $\alpha\dvtx\bN\times\bN\to(0,\infty)$ denotes the
coagulation rate. Formally we have
%
\begin{equation}
\label{eq12}\qquad
\frac{d}{d t}\sum_n\psi(n) f_n =\frac12\sum_{m,n} \alpha
(m,n)f_m(t)f_{n}(t)\bigl(\psi(m+n)-\psi(m)-\psi(n)\bigr)
\end{equation}
for any function $\psi$. An important choice for $\psi$ is $\psi
(n)=n$ with the sum $\sum_n n f_n$ interpreted as
the total mass of particles. For such a choice the right-hand side of
(\ref{eq12}) is $0$ and this is consistent with
our intuition; the total mass for coagulating particles is conserved.
In reality equation (\ref{eq12}) is not valid, and in the case of
$\psi(n)=n$ we only have
%
\begin{equation}
\label{eq13}
\frac{d}{d t}\sum_n n f_n \le0 .
\end{equation}
Analytically speaking, we cannot interchange the differentiation with
the summation in (\ref{eq12}), and such an
interchange can take place only if some suitable restrictions on the
size of the coagulation rate $\alpha(m,n)$ is
imposed as $m$ and $n$ get large. The strict inequality in (\ref
{eq13}) does not contradict the conservation of
mass; for the sufficiently fast growing $\alpha$, particles of infinite
size---the so-called gels---are formed, and
the sum $\sum_n n f_n$ no longer represents the total mass. More
precisely, if we write $g_n=n f_n$
for the total mass of particles of size $n$, then what we really have is
%
\begin{equation}
\label{eq14}
\frac{d}{d t}\Biggl(\sum_{n=1}^{\infty} g_n +g_\infty\Biggr)=0.
\end{equation}

A Marcus--Lushnikov process (MLP) is formulated as a simple microscopic
model to study coagulation
and gelation phenomena.
MLP is a Markov process which is defined on a finite state space $E_N$
given by
\[
E_N=\biggl\{\bL=(L_1,L_2,\ldots,L_n,\ldots)\dvtx \sum_n n L_n=N, 0\le L_n\in
\bZ{\mbox{ for each }} n\biggr\}.
\]
What we have in mind is that $L_n$ is the total number of particles of
size $n$, and the condition $\sum_n n L_n=N$
means that $N$ is indeed\vspace*{1pt} the total mass of particles.
The process $(\bL^{(N)}(t)=\bL(t)\dvtx t\in[0,\infty))$ is a Markov process
with infinitesimal generator $\cA=\sum_{m,n=1}^\infty\cA_{m,n}$,
where
\[
\cA_{m,n}F(\bL)=\frac1{2N}{\alpha(m,n)}\bigl(L_m L_n-\mathbh
{1}(m=n)L_m\bigr)\bigl(F(\bL^{m,n})-F(\bL)\bigr).
\]
When $m\neq n$, $\bL^{m,n}$ is obtained from $\bL=(L_1,L_2,\ldots)$
by replacing $L_n, L_m$ and $L_{n+m}$
with $L_n-1, L_m-1$ and $L_{n+m}+1$, respectively; when $m=n$, $\bL
^{m,n}$ is obtained from $\bL=(L_1,L_2,\ldots)$ by replacing $L_n$ and $L_{2n}$
with $L_n-2$ and $L_{2n}+1$, respectively.\vadjust{\goodbreak}
In words, with rate $\alpha(m,n)/N$, a pair of particles of sizes $m$ and
$n$ is replaced with a single particle of size $m+n$.
Note that the number of such pairs is $L_m L_n$ if $n\neq m$, and this
number becomes $L_n(L_n-1)$ if $m=n$.
Also note that we intentionally have chosen a coagulation rate
proportional to $N^{-1}$. The reason for this has to do with the fact
that all pairs of particles are prone to coagulate, and, as a result,
a typical particle undergoes a huge number of coagulations in one unit
of time as $N$ gets large. Our rescaling
of $\alpha$ guarantees that, on average, a single particle experiences
only a finite number of coagulations. The probability measure
and the expectation associated with the Markov process $\bL(t)$ are
denoted by $\bP_N$ and $\bE_N$, respectively.

The connection between MLP and SE is that the large $N$ limit
$f_n:=\lim_N L^{(N)}_n/ N$ is expected to exist and satisfy SE. For
this, however,
suitable assumptions on $\alpha$ are needed. Before stating these
conditions and a precise theorem relating MLP to SE, let
us make some preparations. Set
\[
E=\biggl\{{\bff}=(f_1,f_2,\ldots,f_n,\ldots)\dvtx \sum_n n
f_n\le1,f_n\ge0 {\mbox{ for each }} n\biggr\}\subset E'=[0,\infty)^\bN.
\]
We equip $E'$ with the product topology. Evidently, $E$ is a compact
subset of~$E'$.
Let us write $\cD=\cD([0,T];E)$ for the Skorohod space of functions
from the interval $[0,T]$ into $E$.
The space $\cD$ is equipped with Skorohod topology. The Markov process
$(\bL(t)\dvtx t\in[0,T])$ induces a probability measure
$\cP_N$ on $\cD$ via the transformation $\bL\mapsto{\bff}$, where
${\bff}=(f_n\dvtx n\in\bN)$, with $f_n=L_n/N$.
We are now ready to state our first result.
%
\begin{theorem}
\label{th11} Assume
%
\begin{equation}\label{eq15}
\sup_{n,m}\frac{\alpha(m,n)}{m+n}<\infty,
\end{equation}
and that initially
%
\begin{equation}\label{eq16}\quad
\lim_{k\to\infty}\limsup_{N\to\infty}\bE_N\frac1N\sum_{n\ge
k}nL_n(0)=0,\qquad
\lim_{N\to\infty}\bE_N \biggl|\frac{L_n(0)}{N}-f^0_n\biggr|=0.
\end{equation}
Then the sequence of probability measures $\{\cP_N\}$ is tight, and if
$\cP$ is a limit point of $\{\cP_N\}$,
then $\cP$ is concentrated on the unique solution to SE subject to the
initial condition ${\bff}(0)={\bff}^0$.
\end{theorem}
%
\begin{remark}\label{Remark11}
The existence of a unique solution to SE under (\ref{eq15}) has been
established in Ball and Carr \cite{BC}. Even though we have not been
able to find a proof of Theorem~\ref{th11} in the literature, we skip
the proof because a straightforward adaption of \cite{BC} can be used
to prove Theorem~\ref{th11}.
\end{remark}

We now turn to the question of gelation which is the primary purpose of
this article. We first recall a result of
Escobedo at al. \cite{EMP} on solutions to SE. We set $M(t)=M({\bff}
,t)=\sum_n n f_n(t)$.\vadjust{\goodbreak}
%
\begin{theorem}
\label{th12}
Assume that $\alpha(m,n)\ge(mn)^{a}$, for some $a>\frac12$. Then there
exists a constant
$C_0(a)$ such that for any solution ${\bff}$ of SE,
%
\begin{equation}\label{eq19}
\int_0^\infty M(t)^2\,dt\le C_0(a)M(0).
\end{equation}
In particular, gelation occurs sometime before $T_{0}=C_0(a)/M(0)$.
That is, for \mbox{$t> T_{0}$}, we have $M(t)<M(0)$.
\end{theorem}

We now discuss the microscopic analog of Theorem~\ref{th12}
for MLP. For this,
let us define stopping times
%
\begin{equation}\label{eq110}
\tau^{(N)}(b,c,\delta)=\tau(b,c,\delta)=\inf\biggl\{t\dvtx N^{-1}\sum_{n\ge
cN^b}nL_{n}(t)\ge\delta\biggr\}.
\end{equation}
The following was established by Jeon \cite{J1}.
%
\begin{theorem}
\label{th13} Assume that $\alpha(m,n)\ge(mn)^{a}$, for some $a>\frac12$.
Then for every $b$ and $\delta\in(0,1)$ and $c>0$,
%
\begin{equation}\label{eq111}
\sup_N\bE_N \tau(b, c,\delta)<\infty.
\end{equation}
\end{theorem}
%
\begin{remark}\label{Remark12}
(i) Section~\ref{sec2} is devoted to the proof of Theorem
\ref{th13}.
Even though we are not introducing any new idea and employing the same
approach as in \cite{J1}, our proof is shorter, more straightforward
and simpler.

(ii) A weaker form of Theorem~\ref{th13} was established by Aldous \cite
{A2} for a special class of coagulation rates $\alpha$.
\end{remark}

Note that if the assumption of Theorem~\ref{th13} holds, then condition (\ref
{eq15}) is no longer true, and, in fact,
we need to modify SE if the sol--gel interaction is significant. It
turns out that if
%
\begin{equation}\label{eq112}
\lim_{m\to\infty}\frac{\alpha(m,n)}{m}=:\bar\alpha(n)
\end{equation}
exists for every $n$, then it is not hard to figure out what the
corrected SE looks like. Under (\ref{eq112}), we still have (\ref
{eq11}), but now
with a modified loss term. More precisely, $Q_n = Q_n^{+} - \hat
Q_n^{-}$, where the modified loss term $\hat Q_n^{-}$ reads as
%
\begin{equation}\label{eq113}
\hat Q_n^{-}(f)(t) = \sum_{m=1}^{\infty} \beta(m,n)g_m(t)g_n(t)+\beta
(n,\infty
)g_n(t)g_\infty(t)
\end{equation}
with $g_n=nf_n$, $\beta(n,m)=\alpha(n,m)/(mn)$, and $\beta(n,\infty)$ measures
the amount of coagulation between particles of size $n$ and gels. When
the condition of Theorems~\ref{th12} or~\ref{th13} holds,
we have that $g_\infty(t)>0$ for $t>T_{\mathrm{gel}}$. In fact, if (\ref{eq112})
holds, then $\beta(n,\infty)$ is simply given by
%
\begin{equation}
\label{eq114}
\beta(n,\infty)= \frac{\bar\alpha(n)}{n}.\vadjust{\goodbreak}
\end{equation}

The analog of Theorem~\ref{th11} in this case is Theorem~\ref{th14}.
%
\begin{theorem}
\label{th14} Assume (\ref{eq112}).
Then the sequence of probability measures $\{\cP_N\}$ is tight.
Moreover, if $\cP$ is a limit point of $\{\cP_N\}$,
then $\cP$ is concentrated on the space of solutions to the modified
SE with the loss term given by
(\ref{eq113}) and (\ref{eq114}) and $g_\infty=1-\sum_ng_n$.
\end{theorem}
%
\begin{remark}\label{Remark13}
(i) Theorem~\ref{th11} under the stronger condition $\bar\alpha(n)=\beta
(n,\infty)=0$
was established in \cite{J1}.
This condition does not exclude gelation. However, even though
a fraction of the density comes from gels (i.e., $g_\infty>0$) after the
gelation time,
the sol--gel interaction is sufficiently weak that can be ignored
in the macroscopic description of the model.

\mbox{}\hphantom{i}(ii) The continuous analog of ML model has been studied in Norris
\cite{N} and Fournier--Giet
\cite{FoG}. In the continuous variant of ML the cluster sizes take
values in $(0,\infty)$ and all $m$
summations in SE (\ref{eq11}), and modified SE are replaced with $dm$
integrations. In
the continuous case, Theorem~\ref{th11} under the stronger condition $\bar
\alpha(n)=0$
was established
in \cite{N} and under the assumption (\ref{eq112}) in \cite{FoG}. As
is stated in
\cite{FoG}, the modified SE has already been predicted by Flory \cite{Fl}.
See also Fournier and Laurencot
\cite{FoL} where a variant of continuous ML with cutoff has been studied.

(iii) It is not hard to understand why a condition like (\ref
{eq112}) facilitates the derivation of the modified
Smoluchowski's equation. The main idea is that even though the function
$\bff\mapsto\sum_{m}\alpha(m,n)f_m$ is not a continuous function with
respect to the product topology whenever $\bar\alpha(n)\neq0$, the
function $\bff\mapsto\sum_{m}(\alpha(m,n)-m\bar\alpha(n))f_m$ is
continuous. This can be easily used to establish Theorem~\ref{th14} by
standard arguments, providing a rather more direct proof of Theorem
\ref{th14}
than the one appeared in \cite{FoG}.

\mbox{}\hspace*{1pt}(iv) If the condition (\ref{eq112}) fails and instead we have
the weaker property
\[
\sup_m \alpha(m,n)/m<\infty,
\]
it is not clear what macroscopic equation, if any describes the
evolution of densities.
\end{remark}

We next address the question of instantaneous gelation. We first recall
a result of Carr and da Costa \cite{CdC}.
%
\begin{theorem}
\label{th15} Assume that for some $q>1$, we have that $\alpha(m,n)\ge
m^q+n^q$. Then $M(t)<M(0)$
for every solution of SE and every $t>0$. In words, gelation occurs
instantaneously.
\end{theorem}

We now state a theorem that is the microscopic analog of Theorem~\ref{th15}.
To this end, let us define
\[
T_k(\delta)=\inf\biggl\{t\dvtx N^{-1}\sum_{n\ge
k}nL_n(t)\ge\delta\biggr\},\quad\ \
\hat T_A^{(N)}(\delta)=\hat T_A(\delta)=T_{A\log N/\log\log
N}(\delta).\nonumber
\]

\begin{theorem}
\label{th16} Assume that $\alpha(m,n)\ge m^q+n^q$, for some $q\in(1,2)$.
Then for every positive $\delta<1$, $A<q(2-q)^{-1}(6-q)^{-1}$ and
$\theta<\bar\eta$,
there exists a constant
$C_2=C_2(q,\theta,A)$,
such that
%
\begin{equation}\label{eq116}
\bE_N \hat T_A(\delta)\le C_2(1-\delta)^{-1}( {\log N})^{-\theta}.
\end{equation}
Here $\bar\eta=\bar\eta(q,A)=\min((q-1)/4,\bar s+q-2)$ with $\bar
s$ given by (\ref{eq35}) below.
\end{theorem}
%
\begin{remark}\label{Remark14}
(i) Note that the condition of Theorem~\ref{th16} is stronger than what
we assume in Theorem~\ref{th13} because
$m^q+n^q\ge2 (mn)^{q/2}$.

\mbox{}\hphantom{i}(ii) Theorem~\ref{th16} is more satisfactory than Theorem~\ref{th15} for three
reasons. On one hand
in Theorem~\ref{th15} we only claim that
if there exists a solution to SE, then such a solution experiences an
instantaneous gelation.
In other words, we are only showing that there is no mass-conserving
solution; however, it is not known if,
under the assumption of Theorem~\ref{th15}, a solution exists.
On the other hand, the macroscopic densities coming from MLP cannot
satisfy (\ref{eq11}) and (\ref{eq113})
because $\beta(n,\infty)=\infty$, and presumably a suitable
modification of SE
would be necessary.
Finally, in Theorem~\ref{th16} we are giving a bound on the time of the
formation of a large particle.
That is, we are giving more information about how instantaneous the
gelation is.
We should mention though that our Lemma~\ref{lem32} in Section~\ref{sec3} is partly
inspired by the proof of Carr and da Costa in \cite{CdC}.

(iii) We note that under the assumption of Theorem~\ref{th13}, the
quickest way for gelation is to wait
first for the creation of several large particles, and then large
particles coagulate among themselves to produce even larger particles
very quickly.
After all if both $m$ and $n$ are of order $\ell$, then $\alpha$ is at
least of order
$\ell^{2a}$ with $2a>1$. However, under the assumption of Theorem~\ref{th16},
gelation is the result of the coagulations
of a large particle with any other particle. Note that for a particle
of size $\ell$ to coagulate with another particle,
it takes a short time of order $O(\ell^{-q})$, and $\sum_{\ell>\ell
_0}\ell^{-q}$ is small if $\ell_0$ is large.
This explains why in Theorem~\ref{th16} we have instantaneous gelation; once a
single large particle is formed, this large
particle coagulates almost immediately with the others to grow even larger.

\mbox{}\hspace*{1pt}(iv) For instantaneous gelation, we only need $\alpha(n,m)\ge\eta
(m)+\eta(n)$ with $\eta$ satisfying
$\sum_n\eta(n)^{-1}<\infty$. A similar comment applies to Theorem
\ref{th17} below.

\mbox{}\hphantom{i}\hspace*{1pt}(v) For simplicity, we avoided the case $q\ge2$. In fact when
$q=2$, (\ref{eq116}) is valid with no restriction on $A$ and $\bar
\eta=1/2$; see Remark~\ref{Remark31} in Section~\ref{sec3}.
The condition $q>2$ leads to instantaneous complete gelation that
will be discussed in Theorem~\ref{th17} below.
\end{remark}

We finally turn to the question of complete gelation. Define
\[
\tilde\tau^{(N)}=\tilde\tau=\inf\{t\dvtx L_{N}(t)=1\}.
\]

\begin{theorem}
\label{th17} Assume that $\alpha(m,n)\ge m^qn+n^qm$, for some $q>1$.
Then there exists a constant $C_3=C_3(q)$ such that
%
\begin{equation}\label{eq119}
\bE_N \tilde\tau\le C_3 \biggl(\frac{\log\log N}{\log N}\biggr)^{q-1}.
\end{equation}
\end{theorem}
%
\begin{remark}\label{Remark15}
In Jeon \cite{J2} it has been shown that a complete instantaneous
gelation occurs if the requirement of Theorem~\ref{th17} is satisfied.
No bound on the time of complete gelation is provided in \cite{J2},
and we believe that our proof is simpler.
\end{remark}

Even though our assumption on $\alpha$ as it appears in Theorem~\ref{th13} is the
most commonly used
condition to guarantee gelation, we now argue that it is the assumption
of Theorem~\ref{th16} that is more physically relevant.
In a more realistic model for the coagulation phenomenon we
would allow spatial dependence for particles. We are now interested on
the evolution of particle density
${\bff}(x,t)=(f_n(x,t)\dvtx n\in\bN)$ where $x\in\bR^d$ represents the
spatial position.
The homogeneous SE is now replaced with the inhomogeneous SE,
\[
\frac{\p}{\p t} f_n(x,t)=\frac12 d(n)\Delta_x f_n(x,t) + Q_n(f)(x,t),
\]
where $d(n)$ denotes the diffusion coefficient of particles of size
$n$, the operator $\Delta_x$ denotes the Laplace operator in $x$
variable and $Q(f)$ has the same form as in the homogeneous SE.
Microscopically, particles have positions, masses and radii. Each particle
travels as a Brownian motion with diffusion
coefficient $d(m)$ where $m$ denotes the mass of the particle.
Particles may coagulate only when they are sufficiently close.
For example, the coagulation occurs between particles of positions $x$
and $x'$ only when $|x-x'|$ is of order
$\e(r+r')$ where $r$ and $r'$ are the radii of particles, and $\e$ is
a small parameter.
When the dimension $d$ is $3$ or more, the initial number of particles
is of order $O(N)$ with $N=\e^{2-d}$.
When particles are close, they coagulate randomly with a rate that is
proportional to $\alpha(m,n)$. This microscopic coagulation rate $\alpha
$ is
not the macroscopic coagulation rate that appears in SE. One can
calculate the macroscopic coagulation rate $\hat\alpha$ from the
microscopic coagulation rate $\alpha$ and the diffusion coefficient
$d(\cdot)$ after
some potential theory. We refer the reader to \cite{HR1,HR2} and \cite
{R} for more details on this model and a precise formula of $\hat\alpha$.
In this model of coagulating Brownian particles, a large microscopic
coagulation rate would not lead to gelation. Instead, the radii of
particles are what matter when it comes to the issue of gelation.
Indeed, if the relationship between the mass $m$ of a particle
and its radius $r$ is given by $r=m^{\chi}$, then for a gelation we
need a condition of the form
$\chi>(d-2)^{-1}$. This is quite understandable in view of Theorem~\ref{th15}
because for a uniformly positive $\alpha$,
the macroscopic coagulation rate $\hat\alpha(m,n)$ behaves like
$(d(m)+d(m))(m^\chi+n^\chi)^{2-d}$ as $m$ and $n$
get large; see~\cite{R}. As a result, if the diffusion coefficients
$(d(n), n\in\bN)$ are uniformly positive\vadjust{\goodbreak} and $\chi>(d-2)^{-1}$,
then $\hat\alpha$ has a super-linear growth as the size of particles get
large. Based on this we conjecture that an
instantaneous gelation would occur if $\chi>(d-2)^{-1}$.

We end this Introduction with the outline of the paper: Section~\ref{sec2} is
devoted to the proof
of Theorem~\ref{th13}. Theorem~\ref{th16} will be established in Section~\ref{sec3}.
Section~\ref{sec4} is devoted to the proof of Theorem~\ref{th17}.

\section{Simple gelation}
\label{sec2}

\mbox{}

\begin{pf*}{Proof of Theorem~\ref{th13}}
For (\ref{eq111}). it suffices to show
that for every $b\in((2a)^{-1},1)$ and positive $\delta$, there exist
constants $C_0(a,b,\delta)$ and $C_0'(a,b,\delta)$ such that
%
\begin{equation}\label{eq21}
\sup_N\bE_N \tau(b, C_0(a,b,\delta),\delta)\le C_0'(a,b,\delta).
\end{equation}
Explicit expressions for the constants $C_0$ and $C'_0$ are given in
(\ref{eq25}) below.

Pick $\beta>0$, and set $\delta_i=\delta+c2^{-i\beta}$ with the constant
$c\in(0,1-\delta]$
so that we always have $\delta_i\le1$. Define the stopping time
\[
T_k=\inf\biggl\{t\dvtx \sum_{n\ge2^i}n L_n(t)\ge\delta_i N \mbox{ for }
i=0,1,\ldots,k\biggr\}
\]
for each $k\in\bN$. Evidently $T_k\le T_{k+1}$. We also define
\[
F_k(\bL)=\frac1N \sum_{n\ge2^{k+1}}n L_n.
\]
By the strong Markov property,
%
\begin{equation}
\label{eq22}
\bE_N F_k(\bL(T_{k+1}))=\bE_N F_k(\bL(T_k))+\bE_N \int
_{T_k}^{T_{k+1}}\cA F_k(\bL(t))\,dt.
\end{equation}
Note that if $T_k\le t< T_{k+1}$, then
\begin{eqnarray*}
\sum_{n\ge2^k}n L_n(t)&\ge&\delta_k N ,\qquad \sum_{n\ge2^{k+1}}n
L_n(t)<\delta_{k+1} N,\\
\frac1N\sum_{n= 2^k}^{2^{k+1}-1}n L_n(t)&\ge&\delta_k-\delta_{k+1}.
\end{eqnarray*}
Let us simply write $\bL=\bL(t)$ with $t$ satisfying $T_k\le t<
T_{k+1}$. For such a configuration $\bL$ we have
that $ \cA F_k(\bL)$ equals
\begin{eqnarray*}
&&\frac1{2N^2}\sum_{m,n} \alpha(m,n)L_m\bigl(L_n-\mathbh
{1}(m=n)\bigr)\\[-2pt]
&&\hphantom{\frac1{2N^2}\sum_{m,n}}
{}\times[(m+n)\mathbh{1}(m+n\ge2^{k+1})
-m\mathbh{1}(m\ge2^{k+1})-n\mathbh{1}(n\ge2^{k+1})]\\[-2pt]
&&\qquad\ge\frac1{2N^2}\sum_{m,n} (mn)^aL_m\bigl(L_n-\mathbh
{1}(m=n)\bigr)\\[-2pt]
&&\qquad\hphantom{\ge\frac1{2N^2}\sum_{m,n}}
{}\times[(m+n)\mathbh{1}(m+n\ge2^{k+1})
-m\mathbh{1}(m\ge2^{k+1})-n\mathbh{1}(n\ge2^{k+1})]\\[-2pt]
&&\qquad= \frac1{2N^2}\sum_{m,n}
(mn)^aL_mL_n[(m+n)\mathbh{1}(m+n\ge2^{k+1})\\[-2pt]
&&\hspace*{128pt}{}-m\mathbh{1}(m\ge2^{k+1})-n\mathbh{1}(n\ge2^{k+1})]\\[-2pt]
&&\qquad\quad{} -\frac1{2N^2}\sum_{m} m^{2a}L_m[(2m)\mathbh{1}(2m\ge2^{k+1})
-2m\mathbh{1}(m\ge2^{k+1})] \\[-2pt]
&&\qquad\ge\frac1{2N^2}\sum_{m,n}
(mn)^aL_mL_n(m+n)\mathbh{1}(m+n\ge2^{k+1}>m,n)\\[-2pt]
&&\qquad\quad{}-\frac1{2N^2}\sum_{m} m^{2a}L_m(2m)\mathbh{1}(2m\ge2^{k+1}>m)\\[-2pt]
&&\qquad= \frac1{2N^2}\sum_{m,n} (mn)^aL_mL_n(m+n)\mathbh{1}(m+n\ge2^{k+1}>m,n)
\\[-2pt]
&&\qquad\quad{}-\frac1{N^2}\sum_{m=2^k}^{2^{k+1}-1} m^{2a+1}L_m\\[-2pt]
&&\qquad\ge\frac1{N^2}\Biggl(\sum_{m=2^k}^{2^{k+1}-1}m^{a+1}L_m\Biggr)
\Biggl(\sum_{n=2^k}^{2^{k+1}-1}n^{a}L_n\Biggr)
-\frac1{N^2}\sum_{m=2^k}^{2^{k+1}-1} m^{2a+1}L_m\\[-2pt]
&&\qquad\ge\frac1{N^2}2^{ka}2^{k(a-1)}2^{-(a-1)^{-}}\Biggl(\sum
_{m=2^k}^{2^{k+1}-1}m L_m\Biggr)^2
- \frac1{N^2}2^{2(k+1)a}\Biggl(\sum_{m=2^k}^{2^{k+1}-1}mL_m\Biggr)\\[-2pt]
&&\qquad\ge2^{k(2a-1)}2^{-(a-1)^{-}}(\delta_k-\delta_{k+1})^2-\frac1N
2^{2a+2ak}\\[-2pt]
&&\qquad=c^2(1-2^{-\beta})^2(2^k)^{2a-1-2\beta}2^{-(a-1)^{-}}-\frac1N2^{2a}(2^k)^{2a}.
\end{eqnarray*}
First we want to make sure that the negative term does not cancel the
positive term. For example, we may try to have
\[
\frac{c^2}2 (1-2^{-\beta})^2(2^k)^{2a-1-2\beta}2^{-(a-1)^{-}}\ge
\frac1N2^{2a}(2^k)^{2a}.
\]
For this it is suffices to assume
\[
2^k\le\bigl(c^2(1-2^{-\beta})^22^{-2a-(a-1)^{-}-1}\bigr)^{1/({1+2\beta
})}N^{1/({1+2\beta})}.
\]
For such integer $k$ we use (\ref{eq22}) to deduce
\[
1\ge\bE_N F_k(\bL(T_{k+1}))\ge\frac{c^2}2(1-2^{-\beta
})^2(2^k)^{2a-1-2\beta}2^{-(a-1)^{-}}
\bE_N (T_{k+1}-T_k).
\]
Hence,
\[
\bE_N (T_{k+1}-T_k)\le2c^{-2}(1-2^{-\beta
})^{-2}2^{(a-1)^{-}}(2^k)^{-(2a-1-2\beta)}.
\]
Summing these inequalities over $k$ yields
\begin{eqnarray*}
\bE_N T_{\ell}&\le&2c^{-2}(1-2^{-\beta})^{-2}2^{(a-1)^{-}}\sum
_{k=0}^{\ell-1}(2^k)^{-(2a-1-2\beta)}\\
&\le&2c^{-2}(1-2^{-\beta})^{-2}2^{(a-1)^{-}}\bigl(1-2^{-(2a-1-2\beta)}\bigr)^{-1}
\end{eqnarray*}
provided that $\beta<a-\frac12$ and
%
\begin{equation}
\label{eq23}
2^\ell\le\bigl(c^2(1-2^{-\beta})^22^{-(a-1)^{-}}2^{-2a-1}\bigr)^{1/({1+2\beta
})}N^{1/({1+2\beta})}.
\end{equation}
If $\ell$ is the largest integer for which (\ref{eq23}) holds, then
\begin{eqnarray*}
2^\ell&\ge&2^{-1}\bigl(c^2(1-2^{-\beta})^22^{-(a-1)^{-}}2^{-2a-1}\bigr)^{
1/({1+2\beta})}N^{1/({1+2\beta})}\\
&=&\!:C(c,a,\beta)N^{1/({1+2\beta})}.
\end{eqnarray*}
From this we deduce
%
\begin{equation}\label{eq24}\quad
\bE_N \tau'_{\beta}\le
2c^{-2}(1-2^{-\beta})^{-2}2^{(a-1)^{-}}\bigl(1-2^{-(2a-1-2\beta
)}\bigr)^{-1}=:C'(c,a,\beta),
\end{equation}
where $\tau'_\beta$ is the first time
\[
N^{-1}\sum_{n\ge k}nL_n(t)\ge\delta
\]
with $k=C(c,a,\beta)N^{1/({1+2\beta})}$.
Since $\beta\in(0,a-\frac12)$ is arbitrary, $b=(1+2\beta)^{-1}$ can take
any value in the interval
$((2a)^{-1},1)$. Finally we choose $c=1-\delta$ to derive (\ref
{eq21}) from (\ref{eq24}) with
%
\begin{eqnarray}\label{eq25}
C_0(a,b,\delta)&=&C\bigl(1-\delta,a,(b^{-1}-1)/2\bigr),\nonumber\\[-8pt]\\[-8pt]
C_0'(a,b,\delta)&=&C'\bigl(1-\delta,a,(b^{-1}-1)/2\bigr).\nonumber
\end{eqnarray}
\upqed\end{pf*}

\section{Instantaneous gelation}
\label{sec3}

This section is devoted to the proof of Theorem~\ref{th16}.
The main ingredient for the proof of Theorem~\ref{th16} is Theorem~\ref{th31}.

\begin{theorem}
\label{th31} Assume that $\alpha(m,n)\ge m^q+n^q$, for some $q\in(1,2)$.
There exist positive constants $C_1=C_1(q,s,\eta,\nu)$ and
$k_0=k_0(q,s,\eta)$
such that if $s>2-q$, $\eta\in(0,(q-1)/4)$, $\delta\in(0,1)$, and
$\nu>1$, then
%
\begin{eqnarray}\label{eq31}
\bE_N T_k(\delta)&\le&
4(2-q)^{-1}k^{-s+2-q}+8k^{s(k-1)+3-q}N^{-1}\nonumber\\
&&{}+C_1(1-\delta)^{-1}k^{-\eta}(\log k)^{1-\eta}\\
&&{}+C_1(1-\delta
)^{-1}k^{3-q/2}(\log k)^3 N^{{-q}/({2s(k-1)})}
\nonumber
\end{eqnarray}
for every $k$ satisfying $k>k_0$ and
%
\begin{equation}\label{eq32}
k^{(k-1)s+2}\le N\le e^{k^\nu},\qquad 2k^{-s}\le1.
\end{equation}
\end{theorem}
%
\begin{remark}\label{Remark31}
For simplicity, we avoided the case $q= 2$. In fact when $q=2$, (\ref
{eq31}) is valid if we replace the first term on the right-hand side
with $4k^{-s}\log k$ [see (\ref{eq316}) below].
\end{remark}

We first demonstrate how Theorem~\ref{th31} implies Theorem~\ref{th16}.
\begin{pf*}{Proof of Theorem~\ref{th16}}
Set $k=A\log N/\log\log N$ in Theorem~\ref{th31}. We note that (\ref{eq32})
is satisfied for large $N$ if
$sA<1$.
Let us first look at the second term on the right-hand side of (\ref
{eq31}). In fact
the second term decays like a negative power of $N$ if $sA<1$. This is because
%
\begin{equation}\label{eq33}
k^{s(k-1)+3-q}N^{-1}\le c_0N^{sA-1}(\log N)^{c_1}
\end{equation}
for some constants $c_0$ and $c_1$. To see this, take the logarithm of
both sides to write
\[
sk\log k +(3-q-s)\log k\le\log c_0+sA\log N+c_1\log\log N.
\]
First select $c_1$ large enough so that
\[
(3-q-s)\log k\le(3-q-s)(\log A+\log\log N)\le c_1\log\log N.
\]
Then observe that if $N$ satisfies $\log\log N\ge A$, then
\[
sk\log k\le sk\log\log N=sA\log N.
\]
This completes the proof of (\ref{eq33}) with $c_0=1$, provided that
$N$ satisfies $\log\log N\ge A$.
Finally we adjust the constant $c_0$ to have the inequality (\ref
{eq33}) even when $N$ satisfies $\log\log N< A$.

We now turn to the last term on the right-hand side of (\ref{eq31}).
By taking the logarithm of the last term,
it is not hard to show that for a positive constant $c_2$,
\[
k^{3- q/2}(\log k)^3 N^{{-q}/({2s(k-1)})}\le c_2 (\log
N)^{3- q/2- q/({2sA})}
(\log\log N)^{q/2}.
\]
The right-hand side of (\ref{eq31}) goes to $0$ as $N\to\infty$, if
\[
s>2-q,\qquad sA<\biggl(\frac{q}{6-q}\biggr)\wedge1.
\]
Now (\ref{eq31}) implies
%
\begin{eqnarray}\label{eq34}
\bE_N T_k(\delta)&\le& c_3(1-\delta)^{-1}[(\log N)^{-\eta}+(\log
N)^{-\eta'}+(\log N)^{-\eta''}]\nonumber\\[-8pt]\\[-8pt]
&&{}\times(\log\log N)^{\g}\nonumber
\end{eqnarray}
with
\[
\eta'=s+q-2,\qquad \eta''=\frac q{2sA}+\frac q2-3,\qquad \g=\max(1,s+q-2).
\]
We now try to optimize (\ref{eq34}) over $s$.
By our assumption on $A$, we know that $(2-q)(6-q)<q/A$.
Choose $s=\bar s$, where
\[
\eta'= \bar s+q-2=\frac q{2\bar sA}+ \frac q2-3=\eta''.
\]
Hence,
%
\begin{eqnarray}\label{eq35}
&\displaystyle \bar
s=\bigl(\sqrt{(1+q/2)^2+2q/A}-1-q/2\bigr)/2,&\nonumber\\[-8pt]\\[-8pt]
&\displaystyle (\bar s-2+q)(2\bar s+6-q)=q/A+(2-q)(q-6)>0.&\nonumber
\end{eqnarray}
As a result $\bar s>2-q$, $\eta'=\eta''>0$, and we can easily see
\[
\bar s A=\frac q{2(\bar s+1+q/2)}<\frac q{6-q}\wedge1
\]
is also valid. In summary,
%
\begin{equation}\label{eq36}
\bE_N T_k(\delta)\le3c_3(1-\delta)^{-1}(\log N)^{-(\eta\wedge\eta
')}(\log\log N)^{\g},
\end{equation}
where $\eta'=\bar s+q-2$ with $\bar s$ as in (\ref{eq35}). Finally
observe that
$\eta\wedge\eta'$ in (\ref{eq36}) can be chosen to be any positive
number $\theta<\bar\eta$.
By decreasing $\theta$ a little bit, we can forget about the double
logarithm and deduce (\ref{eq116}).
\end{pf*}

It remains to establish (\ref{eq31}). The main ingredients for the
proof of Theorem~\ref{th31} are Lemmas~\ref{lem31} and~\ref{lem32}.
Before stating these lemmas and explaining that how they imply Theorem~\ref{th31},
let us provide some heuristics. Perhaps the best way to motivate our
strategy is by taking
a solution $\bff$ of $(\ref{eq11})$ and establishing an instantaneous
gelation for it. This is exactly what Carr and da Costa proved in \cite
{CdC}. However, we offer an alternative proof that is flexible enough
to be carried out microscopically.
The bottom line is that we would like to show that very quickly a good
fraction of particles are large.
We may start with the worst case scenario initially, namely when all
particles are of size
$1$. That is, $f_1(0)=1$ and $f_n(0)=0$ for $n>1$. We then use (\ref
{eq12}) to show that if $M_k(t)=\sum_{n\ge k}nf_n(t)$, then
\[
\frac{dM_{k+1}(t)}{dt}\ge k^{q-1}M_k(t)\bigl(1-M_{k+1}(t)\bigr).
\]
(See the proof of Lemma~\ref{lem31} below.) Note that if $\theta(\delta,k)$ is
the first time $M_{k+1}(t)\ge\delta$,
then for $t<\theta(\delta,k)$,
\[
\frac{dM_{k+1}(t)}{dt}\ge k^{q-1}M_k(t)(1-\delta).
\]
The point is that staring from $M_1(t)=1$ and $M_k(0)=0$ for $k>1$, we
can use induction to deduce
%
\begin{equation}\label{heur1}
M_{k+1}(t)\ge(k!)^{q-2}\bigl( (1-\delta)t\bigr)^k:=\bar\delta_{k+1}(t)
\end{equation}
provided that $t<\theta(\delta,k)$.
What we learn from this is that it takes a short time to have $\bar
\delta_k$
fraction of mass constituting of particles of sizes at least~$k$,\vadjust{\goodbreak}
provided that we choose $\delta_k$
positive but super-exponentially small as $k$ gets larger. As we try to
carry out this argument for $\bL$,
we encounter two difficulties: the discrete nature of the ML model
introduces an additional error coming from
coagulations between two particles of the same size (a microscopic
coagulation rate $L_n^2-L_n$ instead of $L_n^2$),
and the noise in the system. However, the inductive nature of the above
argument allows us to handle these difficulties
and establish a variant of (\ref{heur1}) in Lemma~\ref{lem31}.

Lemma~\ref{lem31} gives us a weak lower bound on the total mass of large
particles because $\bar\delta_{k}$ in (\ref{heur1})
is very small for large $k$. To see how such a weak lower bound can be
improved, let us recall that as in \cite{CdC} we may look at moments
$R_p=\sum_nn^pf_n$ and show that in fact
\[
\frac{dR_{p}(t)}{dt}\ge pR_p(t)^{1+\beta}M_1(t)
\]
with $\beta=(q-1)/(p-1)$. If $t$ is before the gelation time, then
$M_1(t)=1$ and we learn that $R_p(t)$ blows up at
a finite time $t_p$ which is very small if $p$ is very large. Because
of the randomness in our $ML$ model, we do not know
how to work out a microscopic variant of \cite{CdC} argument. Instead
we switch to the moments of large
particles $M_{p,\ell}=\sum_{n\ge\ell}n^pf_n$ and observe that now
\[
\frac{dM_{p,\ell+1}(t)}{dt}\ge pM_{p,\ell}(t)^{1+\beta}\bigl(1-M_\ell
(t)\bigr)M_\ell(t)^{-\beta},
\]
and if $t<\theta(\delta,{\ell}-1)$, then
%
\begin{equation}\label{heur2}
\frac{dM_{p,\ell+1}(t)}{dt}\ge pM_{p,\ell}(t)^{1+\beta}(1-\delta
)\delta^{-\beta}.
\end{equation}
The point is that now the right-hand side of (\ref{heur2}) depends on
the previous $M_{p,\ell}$, and therefore an inductive argument can be
used to show that $M_{p,\ell}(t)$ can get very large for a time $t$
that is small and $p$ that is large. In other words, instead of showing
that $R_p$ becomes infinite at a time $t_p$ that is small, we would
rather show that $M_{p,\ell}(t)$ gets extremely large very quickly. The
inductive nature of (\ref{heur2}) makes it very useful in its
microscopic form. More precisely, in the case of ML process we can show
that a variant of (\ref{heur2}) is true for the $\bL$ process provided
that we take the expectation of both sides. Then by induction on $\ell$
we can show that $M_{p,\ell}(t)$ gets very large very quickly. This is
exactly the role of Lemma~\ref{lem32} below. In fact the induction
starts from $\ell=k$, and we use Lemma~\ref{lem31} to argue that
$M_{p,k}(t)$ is already large for some small $t$ provided that $p$ is
sufficiently large. With the aid of Lemma~\ref{lem32}, we show that if we wait
for another short period of time, either a good fraction of particles
are large, or else the high moments of density become
super-exponentially large in k. Then a crude bound on moments of
particle density demonstrates that the second alternative cannot occur,
and hence gels have already been formed.

To prepare for the statement of the first lemma,
we take a sequence $(\delta_\ell\dvtx\ell=1,\ldots,k)$, and define
\[
\s_\ell=\inf\biggl\{t\dvtx\frac1N\sum_{n\ge r}nL_n(t)\ge\delta_r{\mbox{
for }} r=1, 2,\ldots,\ell\biggr\}.
\]

\begin{lemma}
\label{lem31}
For every decreasing sequence $(\delta_\ell\dvtx\ell=1,\ldots,{k})$
which satisfies
%
\begin{equation}\label{eq37}
\delta_1=1,\qquad 2\delta_{2}\le1\quad \mbox{and}\quad \frac{8 {k}}N\le\delta_{{k}},
\end{equation}
we have
%
\begin{equation}
\label{eq38}
\bE_N\s_{{k}}\le4\sum_{\ell=1}^{{k}-1} \ell^{1-q}\frac{\delta
'_{\ell+1}}{\delta_{\ell}},
\end{equation}
where $\delta'_{\ell+1}=\delta_{\ell+1}+2\ell N^{-1}$.
\end{lemma}

Define
\[
T_{p,r}(A)=\inf\biggl\{t\dvtx\frac1N\sum_{n\ge r}n^pL_n(t)\ge A\biggr\}.
\]
Recall that we simply write $T_{r}(A)$ for $T_{p,r}(A)$ when $p=1$.
%
\begin{lemma}
\label{lem32}
Let $\{m_\ell\dvtx k\le\ell\le h\}$ be an increasing sequence, and pick
$p\ge2$, $\delta>0$. Assume that
$Nm_{\ell+1}\ge p\ell^p$
for every $\ell$, and write $\tau_r$ for $T_{p,r}(m_r)\wedge
T_k(\delta)$.
Then for $h>k$,
%
\begin{equation}
\label{eq39}
\bE_N (\tau_h-\tau_k)\le\frac{2}{1-\delta}\sum_{\ell=k}^{h-1}
\biggl[\frac{\delta^\beta m_{\ell+1}}{pm_\ell^{\beta+1}}+\ell^2\biggl(\frac{p\ell
}{Nm_{\ell+1}}\biggr)^{{q}/({p-1})} \biggr],
\end{equation}
where $\beta=(q-1)/(p-1)$ with $q$ as in the statement of Theorem~\ref{th16}.
\end{lemma}
\begin{pf*}{Proof of Lemma~\ref{lem31}}
To bound the stopping time $\s_\ell$, we use the strong Markov
property to write
%
\begin{equation}
\label{eq310}
\bE_N G'_{\ell+1}(\bL(\s_{\ell+1}))=\bE_N G'_{\ell+1}(\bL(\s
_\ell))+\bE_N \int^{\s_{\ell+1}}_{\s_\ell}\cA G'_{\ell+1}(\bL(t))\,dt,\hspace*{-35pt}
\end{equation}
where $G'_{\ell+1}(\bL)=G_{\ell+1}(\bL)\wedge\delta'_{\ell+1}$, with
\[
G_k(\bL)=\frac1N\sum_{n\ge k}nL_n.
\]
Assume that $\s_\ell<\s_{\ell+1}$, and set
\[
\D_{m,n}=(m+n)\mathbh{1}(m+n\ge\ell+1)-n\mathbh{1}(n\ge\ell
+1)-m\mathbh{1}(m\ge\ell+1).\vadjust{\goodbreak}
\]
We certainly have that $\cA G'_{\ell+1}(\bL)$ is bounded below by
\begin{eqnarray*}
&& \frac1{2N^2}\sum_{m,n}(m^q+n^q) L_n\bigl(L_m-\mathbh{1}(m=n)\bigr)
\mathbh{1}\bigl(G_{\ell+1}(\bL)+N^{-1}\D_{m,n}\le\delta'_{\ell+1}\bigr)\D
_{m,n}\\
&&\qquad\ge\frac1{N^2}\sum_{m,n}(m^q+n^q) L_n\bigl(L_m-\mathbh{1}(m=n)\bigr)
\\
&&\hphantom{\frac1{N^2}\sum_{m,n}}
\qquad\quad{}\times\mathbh{1}\bigl(G_{\ell+1}(\bL)+(m+\ell)/N\le\delta'_{\ell+1}\bigr)m
\mathbh{1}(n\ge\ell\ge m) .
\end{eqnarray*}
Here we restricted the summation to the cases $n\ge\ell\ge m$ and
$m\ge\ell\ge n$ and used symmetry to consider the former case only.
We note that if $n\ge\ell\ge m$, then either $\D_{m,n}=m$
or $m+\ell$. Also note that if $\bL=\bL(t)$ for some $t\in(\s_\ell
,\s_{\ell+1})$
and $n\ge\ell\ge m$, then $G_{\ell+1}(\bL)\le\delta_{\ell+1}$
and $G_{\ell+1}(\bL)+(m+\ell)/N\le\delta'_{\ell+1}$. Hence for such
a configuration $\bL$,
\begin{eqnarray*}
\cA G'_{\ell+1}(\bL)&\ge&\frac1{N^2}\sum_{m,n}(m^q+n^q)m\mathbh
{1}(n\ge\ell\ge m)L_mL_n-\frac2{N^2}\ell^{q+1}
L_\ell\\
&\ge&\frac{1}{N^2}\biggl(\sum_{n\ge\ell}n^qL_n\biggr)\biggl(\sum_{m<\ell
+1}mL_m\biggr)-\frac{2\ell^q}{N^2}\sum_m m L_m\\
&\ge&\ell^{q-1}G_\ell(\bL)\bigl(1-G_{\ell+1}(\bL)\bigr)-\frac{2\ell^q}{N}.
\end{eqnarray*}
If $\s_\ell\le t< \s_{\ell+1}$, then $G_{\ell}(\bL(t))\ge\delta
_{\ell}$, and
$1-G_{\ell+1}(\bL(t))\ge1-\delta_{\ell+1}\ge1/2$ for $\ell\ge
1$, because by our assumption (\ref{eq37}),
$\delta_{\ell+1}\le1/2$. Hence
\[
\cA G'_{\ell+1}(\bL)\ge\tfrac12 \ell^{q-1}\delta_{\ell}-2\ell
^{q}N^{-1}\ge\tfrac14 \ell^{q-1}\delta_{\ell},
\]
where we have used assumption (\ref{eq37}) for the second inequality.
From this and (\ref{eq310}) we deduce
\[
\tfrac14 \ell^{q-1}\delta_{\ell} \bE_N(\s_{\ell+1}-\s_{\ell
})\le
\bE_N [G'_{\ell+1}(\bL(\s_{\ell+1}))- G'_{\ell+1}(\bL(\s_\ell
))]\le\delta'_{\ell+1}.
\]
As a result,
\[
\tfrac14 \ell^{q-1}\delta_{\ell} \bE_N(\s_{\ell+1}-\s_{\ell
})\le
\delta'_{\ell+1},
\]
Hence
\[
\bE_N(\s_{\ell+1}-\s_{\ell})\le4\ell^{1-q}\frac{\delta'_{\ell
+1}}{\delta_{\ell}}.
\]
Summing this inequality over $\ell$ and remembering that $\s_1=0$,
leads to (\ref{eq38}).
\end{pf*}
\begin{pf*}{Proof of Lemma~\ref{lem32}}
\textit{Step} 1.
We note that since $m_\ell<m_{\ell+1}$, we have that $\tau_\ell\le
\tau_{\ell+1}$. Fix some positive $n_0\in\bR$, and write $\theta$ for
the first time $L_n\neq0$
for some $n\ge n_0$.\vadjust{\goodbreak} We also set $\tau'_{\ell+1}=\tau_{\ell
+1}\wedge(\theta\vee\tau_\ell)$.
We use the strong Markov property to write
%
\begin{eqnarray}
\label{eq311}
\bE_N M'_{p,\ell+1}(\bL(\tau'_{\ell+1}))&=&\bE_N M'_{p,\ell+1}(\bL
(\tau_\ell))\nonumber\\[-8pt]\\[-8pt]
&&{}+\bE_N \int^{\tau'_{\ell+1}}_{\tau_\ell}\cA
M'_{p,\ell+1}(\bL(t))\,dt,\nonumber
\end{eqnarray}
where
\[
M_{p,r}(\bL)=\frac1N\sum_{n\ge r}n^pL_n,\qquad
M'_{p,r}(\bL)=M_{p,r}(\bL)\wedge(2m_{r}) .
\]
Write
\[
\Delta_{m,n}:=N^{-1}[(m+n)^p-n^p-m^p]\ge N^{-1}pn^{p-1}m=:\Delta'_{m,n}.
\]
(Here we have used our assumption $p\ge2$.)
We certainly have that the expression $\cA M'_{p,\ell+1}(\bL)$ is
bounded below by
\begin{eqnarray*}
&& \frac1{2N}\sum_{m,n}(m^q+n^q)L_n\bigl(L_m-\mathbh{1}(m= n)\bigr)
\mathbh{1}(m\ge\ell>n{\mbox{ or }}n\ge\ell>m)\\
&&\hphantom{\frac1{2N}\sum_{m,n}}
{} \times[(M_{p,\ell+1}+\Delta_{m,n})\wedge(2m_{\ell+1})-
M_{p,\ell+1}\wedge(2m_{\ell+1})]\\
&&\qquad\ge\frac1{N}\sum_{m,n}(m^q+n^q)L_nL_m\mathbh{1}(n\ge\ell>m)\\
&&\qquad\hphantom{\ge\frac1{N}\sum_{m,n}}
{}\times[(M_{p,\ell+1}+\Delta'_{m,n})\wedge(2m_{\ell+1})-
M_{p,\ell+1}\wedge(2m_{\ell+1})]\\
&&\qquad\ge\frac1{N^2}\sum_{m,n}(m^q+n^q)pn^{p-1}mL_mL_n\mathbh{1}(n\ge
\ell> m)\\
&&\qquad\hphantom{\ge\frac1{N^2}\sum_{m,n}}
{}\times\mathbh{1}(M_{p,\ell+1}+\Delta'_{m,n}\le2m_{\ell+1}).
\end{eqnarray*}
We now assume that $ m<\ell$ and that $\bL=\bL(t)$ for some $\tau
_\ell<t< \tau'_{\ell+1}$.
For such $m$ and $\bL$, we
have
\[
M_{p,\ell+1}(\bL)+\Delta'_{m,n}\le m_{\ell+1}+N^{-1}pn_0^{p-1}\ell
\le2m_{\ell+1}
\]
provided that we choose
\[
n_0=\biggl(\frac{Nm_{\ell+1}}{p\ell}\biggr)^{1/({p-1})}.
\]
For such choices of $\bL$ and $n_0$, we deduce
\begin{eqnarray*}
\cA M'_{p,\ell+1}(\bL)&\ge&\frac{p}{N^2}\biggl(\sum_{n\ge\ell
}n^{p+q-1}L_n\biggr)\biggl(\sum_{m<\ell}mL_m\biggr)\\
&\ge &pM_{p+q-1,\ell}(\bL)\bigl(1-G_{\ell}(\bL)\bigr).
\end{eqnarray*}
If $ t<T_k(\delta)$ and $k\le\ell$, then $G_{\ell}(\bL(t))\le
G_{k}(\bL(t))<\delta$, and
$1-G_{\ell}(\bL(t))\ge1-\delta$. Hence
%
\begin{equation}
\label{eq312}
\cA M'_{p,\ell+1}(\bL)\ge p(1-\delta)M_{p+q-1,\ell}(\bL),
\end{equation}
whenever $\bL=\bL(t)$ for some $t<T_k(\delta)$. On the other hand,
by H\"older's inequality,
\[
M_{p+q-1,\ell}(\bL)=\frac1N\sum_{n\ge\ell}n^{p+q-2}nL_n\ge G_\ell
(\bL)^{-\beta}M_{p,\ell}^{1+\beta},
\]
where $\beta=(q-1)/(p-1)$.
From this and (\ref{eq312}) we deduce that if $\bL=\bL(t)$ for some
$t$ satisfying $t\in(\tau_\ell, \tau'_{\ell+1})$, then
\[
\cA M'_{p,\ell+1}(\bL)\ge{ p(1-\delta)}\delta^{-\beta}M_{p,\ell
}^{\beta+1}(\bL)\ge{ p(1-\delta)}\delta^{-\beta}m_{\ell}^{\beta+1}.
\]
Here we have used the fact that if $\tau_\ell< \tau'_{\ell+1}$,
then $\tau_{\ell}=T_{p,\ell}(m_\ell)$. [Simply because if $\tau
_\ell\neq T_{p,\ell}(m_\ell)$, then we must have that $T_k(\delta
)<T_{p,\ell}(m_\ell)$, which implies that
$\tau'_{\ell+1}=\tau_{\ell}=T_k(\delta)$, $ \tau'_{\ell+1}-\tau
_{\ell}=0$.]
This and (\ref{eq311}) imply
\begin{eqnarray*}
{ p(1-\delta)}\delta^{-\beta}m_{\ell}^{\beta+1} \bE_N(\tau'_{\ell
+1}-\tau_{\ell})&\le&
\bE_N [M'_{p,\ell+1}(\bL(\tau'_{\ell+1}))- M'_{p,\ell+1}(\bL
(\tau_\ell))]\\
&\le&2m_{\ell+1}.
\end{eqnarray*}
Therefore,
\[
\bE_N(\tau'_{\ell+1}-\tau_{\ell})\le\frac{2\delta^\beta
}{p(1-\delta)}\frac{m_{\ell+1}}{m_\ell^{\beta+1}}.
\]

Hence for (\ref{eq39}) it suffices to establish
%
\begin{equation}
\label{eq313}
\bE_N (\tau_{\ell+1}-\tau'_{\ell+1})\le\frac{2\ell^2}{1-\delta}
\biggl(\frac{Nm_{\ell+1}}{p\ell}\biggr)^{{-q}/({p-1})} .
\end{equation}

\textit{Step} 2. To establish (\ref{eq313}), observe that if $\tau
_{\ell+1}>\tau'_{\ell+1}$, then
the configuration $\bL(\tau'_{\ell+1})$ has at least one particle of
size $n\ge n_0$.
Let us mark one such particle and follow its interaction with other
particles for $t\ge\tau'_{\ell+1}$.
When this particle coagulates with any other particle of size $a$, then
we increase its size
$n(t)$ by $a$ and remove the other particle from the system. We write
$\beta_1<\beta_2<\cdots$ for the consecutive coagulation times of the
marked particle with
particles of sizes $m<\ell$. Let us define an auxiliary process
$(Z(t),K(t))$ that is defined for $t\ge\tau'_{\ell+1}$ with
$Z(\tau'_{\ell+1})=K(\tau'_{\ell+1})=0$ and each time our marked
particle coagulates with a particle of size $m<\ell$,
the value of $K$ increases by $1$, and the value of $Z$ increases
by $pn_0^{p-1}N^{-1}$. So, the process $K$ simply counts the number of
such coagulations, and
$Z(t)=pn_0^{p-1} N^{-1}K(t)$. Since at such a coagulation, the
expression $M_{p,\ell+1}$ increases by
$\Delta_{m,n(t)}\ge pn_0^{p-1}mN^{-1}\ge pn_0^{p-1}N^{-1}$,
with $n$ denoting the size of the marked particle, we have
\[
M_{p,\ell+1}(\bL(\beta_j))\ge jpn_0^{p-1}N^{-1}= jm_{\ell+1}\ell^{-1}.\vadjust{\goodbreak}
\]
The right-hand side is $m_{\ell+1}$ if $j=\ell$. As a result, $\beta
_\ell\ge T_{p,\ell+1}(m_{\ell+1})$ and (\ref{eq313}) would follow
if we can show
%
\begin{eqnarray}
\label{eq314}
\bE_N (\tau_{\ell+1}-\tau'_{\ell+1})&\le&\bE_N \bigl(\beta_{\ell}\wedge
T_k(\delta)-\tau'_{\ell+1}\bigr)\nonumber\\[-8pt]\\[-8pt]
&\le&\frac{2\ell^2}{1-\delta}
\biggl(\frac{Nm_{\ell+1}}{p\ell}\biggr)^{{-q}/({p-1})} .\nonumber
\end{eqnarray}
For this, use the Markov property to write
\begin{eqnarray*}
\ell&=&\bE_N \bigl(K(\beta_\ell)-K(\tau'_{\ell+1})\bigr)
\ge\bE_N \bigl(K\bigl(\beta_\ell
\wedge T_k(\delta)\bigr)-K(\tau'_{\ell+1})\bigr)\\
&=&\bE_N \int^{\beta_\ell\wedge T_k(\delta)}_{\tau'_{\ell+1}}
\frac{1}{2N}\sum_{m<\ell}\alpha(m,n(t))\bigl(L_m(t)-\mathbh{1}(n(t)=m)\bigr)
\,dt\\
&=&\bE_N \int^{\beta_\ell\wedge T_k(\delta)}_{\tau'_{\ell+1}}
\frac{1}{2N}\sum_{m<\ell}\alpha(m,n(t))L_m(t) \,dt\\
&\ge&\bE_N\int^{\beta_\ell\wedge T_k(\delta)}_{\tau'_{\ell+1}}\frac
{ n_0^q}{2\ell N}\sum_{m<\ell}mL_m(t)\,dt \\
&\ge&\frac{ n_0^q}{2\ell}\int^{\beta_\ell\wedge T_k(\delta)}_{\tau
'_{\ell+1}}\bigl(1-G_\ell(\bL(t))\bigr)\,dt\\
&\ge&\frac{(1-\delta) n_0^q}{2\ell}\bE_N
\bigl({\beta_\ell\wedge T_k(\delta)}-{\tau'_{\ell+1}}\bigr),
\end{eqnarray*}
where $n(t)$ denotes the size of the marked particle. Here the third
equality requires an explanation:
Recall that by our assumption $Nm_{\ell+1}\ge p\ell^p$, which implies that
$n(t)\ge n_0\ge\ell$ and $\mathbh{1}(n(t)=m)=0$ for $m<\ell$.
Hence (\ref{eq314}) is true, and this completes the proof of (\ref{eq39}).
\end{pf*}

We are now ready to prove Theorem~\ref{th31}.
\begin{pf*}{Proof of Theorem~\ref{th31}}
\textit{Step} 1. There are various parameters in Lemmas~\ref{lem31}
and~\ref{lem32} that
we need to choose to serve our goal.
We start from specifying the sequence $\{\delta_{\ell}\dvtx\ell=1,\ldots
,{k}\}$. We are going to choose
$\delta_\ell={{k}}^{-s(\ell-1) }$.
Note that the conditions in (\ref{eq37}) hold if
%
\begin{equation}\label{eq315}
8 {{k}}^{s {{k}}-s+1}\le N,\qquad 2k^{-s}\le1.
\end{equation}

By (\ref{eq38}),
%
\begin{eqnarray}\label{eq316}
\bE_N \s_{{k}} &\le&4 {{k}}^{-s}\sum_{\ell=1}^{{{k}}-1}( \ell
^{1-q}+2N^{-1}\ell^{2-q}k^{s\ell})\nonumber\\[-8pt]\\[-8pt]
&\le&4(2-q)^{-1} {{k}}^{-s+2-q}+8k^{s(k-2)+3-q}N^{-1},\nonumber
\end{eqnarray}
because $ {\delta_{\ell+1}}/{\delta_{\ell}}= {{k}}^{-s}$.\vadjust{\goodbreak}

\textit{Step} 2.
We now would like to apply Lemma~\ref{lem32}. For this we first specify $p$ to
be $2s({k}-1)+1$. We note that $p>2$
because $k>1+(2s)^{-1}$ follows from the condition $2k^{-s}\le1$ of
(\ref{eq315}).
Also note that if $\bL=\bL(\s_{{k}})$, then
%
\begin{equation}\label{eq317}
M_{p,{k}}(\bL)= \frac1N\sum_{n\ge{{k}}}n^pL_n\ge{k}^{p-1}\delta
_{{{k}}}={{k}}^{p-1-s({k}-1)}
={{k}}^{s({k}-1)}.
\end{equation}
Because of this, we are going to set $m_{{k}}={{k}}^{s({k}-1)}$,
so that $T_{p,{k}}(m_{{k}})\le\s_{{k}}$.
We then specify $m_\ell$ for $\ell>{k}$. We require that $m_{\ell
+1}m_{\ell}^{-\beta-1}=\ell^{-\eta}$
for some $\eta\in(0,1)$. This requirement leads to the formula
%
\begin{equation}\label{eq318}
m_{\ell}=m_{{k}}^{(\beta+1)^{\ell-{k}}}\prod_{r={k}}^{\ell-1}
r^{-\eta(\beta+1)^{\ell-r-1}}.
\end{equation}
In order to apply Lemma~\ref{lem32}, we need to check that $\{m_\ell\dvtx k\le\ell
\le h\}$ is an increasing sequence
and that $Nm_{\ell+1}\ge p\ell^p$.
We establish this assuming that $h=Ak\log k$, and $k$ is sufficiently
large. Since $m_{\ell+1}/m_\ell=m_\ell^\beta\ell^{-\eta}$,
we only need to show that $m_\ell^\beta>\ell^\eta$ for the
monotonicity of $m_\ell$.
Note that for $m_k^\beta>k^\eta$, we need to assume that
$\eta<(q-1)/2$. As we will see shortly, for $m_\ell^\beta>\ell^\eta$
for $k\le\ell\le h$, with $h=Ak\log k$
we need to assume more; it suffices to have $\eta<(q-1)/4$.

Observe
\[
\log m_{\ell}={(\beta+1)^{\ell-{k}}}\Biggl[\log m_{{k}}-\eta\sum
_{r={k}}^{\ell-1}(\beta+1)^{k-r-1}
\log r\Biggr].
\]
Let us write $a=\log(\beta+1)$. Note that for sufficiently large
$k>k_1((q-1)/s)$,
the function $r\mapsto e^{-(r+1)a}\log r$ is decreasing over the
interval $[k,\infty)$.
As a result,
%
\begin{eqnarray}\label{eq319}
\sum_{r={k}}^{\ell-1}(\beta+1)^{k-r-1}
\log r &=&\sum_{r={k}}^{\ell-1}e^{(k-r-1)a}\log r\le e^{ka}\int
_{{k}}^\infty e^{-r a}
\log( r-1)\,dr\nonumber\\
&\le& e^{ka}\int_{{k}}^\infty e^{-r a}
\log r \,dr\nonumber\\
&=&a^{-1}\log{k}+a^{-1} e^{ka}\int_{{k}}^\infty e^{-r a}
r^{-1}\,dr\\
&=& a^{-1}\log{k}+ a^{-1}e^{ka}\int_{a{k}}^\infty e^{-r }
r^{-1}\,dr\nonumber\\
&\le& a^{-1}\log{k}+a^{-1}e^{ka}\bigl(\log^+ (a{k})^{-1}+c_2 \bigr),\nonumber
\end{eqnarray}
where $c_2=\int_1^\infty e^{-r}r^{-1}\,dr$. Here we integrated by parts for
the second equality. Recall that $a=\log(\beta+1)$ with $\beta
=(q-1)/(p-1)$ and $p-1=2s(k-1)$. As a result, $a{k}$ is bounded and
bounded away from $0$, and
\[
\lim_{{k}\to\infty}\frac{a^{-1}\log{k}+ a^{-1}e^{ka}(\log^+ (a{k})^{-1}+c_2)}
{{k}\log{k}}=\frac{2s}{q-1}.
\]
From all this we learn
%
\begin{equation}\label{eq320}
\liminf_{{k}\to\infty}{(\beta+1)^{-(\ell-{k})}}\frac{\log m_{\ell
}}{{k}\log{k}}\ge s-\frac{2\eta s}{q-1}.
\end{equation}
We choose $\eta\in(0,(q-1)/2)$ so that the left-hand side is
positive. For such~$\eta$, choose $\g$ such that
\[
s\biggl(1-\frac{2\eta}{q-1}\biggr)>\g>0.
\]
Hence, for sufficiently large ${k}>k_2(q,s,\g)$ and every $\ell>k$,
%
\begin{equation}\label{eq321}
\log m_\ell\ge\g{(\beta+1)^{\ell-{k}}}{k}\log{{k}},
\end{equation}
which implies
%
\begin{equation}\label{eq322}
\beta\log m_\ell\ge\frac{\g(q-1)}{2s}{(\beta+1)^{\ell-{k}}}\log{{k}}.
\end{equation}
Note that $k_2$ is independent of $\ell$ because (\ref{eq320})
follows from (\ref{eq319})
and the right-hand side of (\ref{eq319}) is independent of $\ell$.
For the monotonicity of the sequence $\{m_\ell\dvtx k\le\ell\le h\}$, we
need to show that $\beta\log m_\ell>\eta\log\ell$.
By (\ref{eq322}), it suffices to have
%
\begin{equation}\label{eq300}
\beta\log m_\ell\ge\frac{\g(q-1)}{2s}{(\beta+1)^{\ell-{k}}}\log
{{k}}>\eta\log h\ge\eta\log\ell.
\end{equation}
Since $h=Ak\log k$, it suffices to have
\[
\frac{\g(q-1)}{2s}{(\beta+1)^{\ell}}\log{{k}}>\eta{(\beta
+1)^{{k}}}[\log k+\log\log k +\log A]
\]
for $\ell\ge k$. This is true if $k>k_3(q,s,\g,\eta,A)$ for a
suitable $k_3$ and $\g(q-1)/\break(2s)>\eta$. As a result, we need to
select $\g$ such that
%
\begin{equation}\label{eq323}
s\biggl(1-\frac{2\eta}{q-1}\biggr)>\g>\frac{2s\eta}{q-1}.
\end{equation}
Such a number $\g$ exists if $\eta\in(0,(q-1)/4)$. So, let us assume
that $\eta\in(0,(q-1)/4)$ and choose $\g=s/2$.
In summary, there exists a constant $k_4=k_4(q,s,\eta,A)$ such that
if $k>k_4$ and $h=Ak\log k$, then the sequence $(m_\ell\dvtx\break\ell=k,\ldots
,h)$ is increasing.

\textit{Step} 3. So far we know that $m_\ell$ is increasing. In order
to apply Lemma~\ref{lem32}, we still need to check that
$Nm_{\ell+1}\ge p\ell^p$ for $\ell$ satisfying $k\le\ell\le h$.
We establish this by induction on $\ell$. If $\ell=k$, then what we
need is
\begin{eqnarray*}
Nm_{k+1}&=&Nm_k m_k^\beta k^{-\eta}= Nk^{s(k-1)} k^{(q-1)/2}k^{-\eta} \ge
pk^p\\
&=&\bigl(2s(k-1)+1\bigr)k^{2s(k-1)+1}.
\end{eqnarray*}
Since $(q-1)/2>\eta$, it suffices to have
%
\begin{equation}\label{eq30}
N\ge k^{s(k-1)+2}
\end{equation}
and $k\ge k_5(q,s,\eta)$.

We now assume that $Nm_{\ell}\ge p(\ell-1)^p$ is valid and try to deduce
$Nm_{\ell+1}\ge p\ell^p$. Indeed
\[
Nm_{\ell+1}=Nm_\ell m_\ell^\beta\ell^{-\eta}\ge p(\ell-1)^p m_\ell
^\beta\ell^{-\eta}
\]
by induction hypothesis, and this is greater than $p\ell^p$, if
\[
m_\ell^\beta\ge\ell^{\eta}\biggl(1+\frac1{\ell-1}\biggr)^p
\quad\mbox{or}\quad \beta\log m_\ell\ge\eta\log\ell+p\log\biggl(1+\frac1{\ell
-1}\biggr) .
\]
Since $p=2s(k-1)+1$, the second term on the right-hand side is bounded,
and we only need to verify
%
\begin{equation}\label{eq3000}
\beta\log m_\ell\ge\eta\log h+c_3\ge\eta\log\ell+c_3
\end{equation}
for a constant $c_3$. Except for the extra constant $c_3$, this is
identical to (\ref{eq300}) and we can readily see that
condition (\ref{eq323}) would guarantee (\ref{eq3000}) if $k\ge
k_6(q,s,\g,\eta,A)$. In summary,
$Nm_{\ell+1}\ge p\ell^p$ is valid for $\ell$ satisfying $k\le\ell
\le h$, provided that $k$ is sufficiently large,
and (\ref{eq30}) is satisfied. We observe that (\ref{eq30}) implies
the first inequality in (\ref{eq315})
for $k\ge8$.

\textit{Step} 4. We assume that $\eta\in(0,(q-1)/4)$ and that $\g=s/2$.
As before, we set $\tau_r=T_{p,r}(m_r)\wedge T_k(\delta)$. Since
$\tau_{{k}}\le T_{p,{k}}(m_{{k}})\le\s_{{k}}$, we may apply Lemma~\ref{lem32} to assert,
\begin{eqnarray*}
\bE_N (\tau_{h}-\s_{{k}}) &\le& \frac{2\delta^\beta}{p(1-\delta
)}\sum_{\ell={k}}^{h-1} \ell^{-\eta}+
\frac{2N^{{-q}/({p-1})}}{(1-\delta)}
\sum_{\ell={k}}^{h-1}\ell^2\biggl(\frac{p\ell}{m_{\ell+1}}\biggr)^{
{q}/({p-1})}\\
&\le&\frac{2\delta^\beta}{p(1-\delta)(1-\eta)}h^{1-\eta}+\frac
{2(m_{k}N)^{{-q}/({p-1})}}{(1-\delta)}
\sum_{\ell={k}}^{h-1}\ell^2({p\ell})^{{q}/({p-1})}\\
&\le&\frac{2}{p(1-\delta)(1-\eta)}h^{1-\eta}+\frac
{2(m_{k}N)^{{-q}/({p-1})}}{(1-\delta)}
h^3(ph)^{{q}/({p-1})}.
\end{eqnarray*}
Hence
%
\begin{eqnarray}\label{eq324}
\bE_N (\tau_{h}-\s_{{k}}) &\le&2\bigl({p(1-\delta)(1-\eta
)}\bigr)^{-1}h^{1-\eta}\nonumber\\[-8pt]\\[-8pt]
&&{}+2(1-\delta)^{-1}k^{-q/2}N^{{-q}/({2s(k-1)})}
h^3({ph})^{{q}/({p-1})}.\nonumber
\end{eqnarray}

Our strategy is to choose $h$ sufficiently large so that $\tau
_{h}=T_{k}(\delta)$, because we are interested in bounding $T_k(\delta
)$. We have the trivial bound $M_{p,h}\le N^{p-1}$ because $N^{-1}\sum
_n nL_n=1$. Hence if $h$ is sufficiently large so that $m_{h}>N^{p-1}$,
then $T_{p,h}(m_h)=\infty$ and
as a result $\tau_{h}=T_{k}(\delta)$. For $m_{h}>N^{p-1}$,
we need
%
\begin{equation}\label{eq325}
\log m_{h}>(p-1)\log N=2s({k}-1)\log N.\vadjust{\goodbreak}
\end{equation}
By (\ref{eq321}),
\[
\log m_h\ge\g{(\beta+1)^{h-{k}}}{k}\log{{k}}
\]
for ${k}>k_2(q,s,\g)$ and $\g=s/2$. As a result, the condition $m_{h}>N^{p-1}$
is realized if
\[
\g{(\beta+1)^{h-{k}}}{k}\log{{k}}\ge2s{k}\log N
\]
or equivalently
\[
\log\g+Ak\log(\beta+1)\log k-k\log(\beta+1)+\log\log{k}\ge\log(2
s)+\log\log N.
\]
Since $\lim_k k\log(\beta+1)=(q-1)/(2s)$, as $k$ goes to infinity, we
pick $\mu\in(0,(q-1)/{2s})$ and choose $k_7((q-1)/s)$ so that if
$k>k_7$, then $k\log(\beta+1)>\mu$. For such~$k$, we only need to have
\[
\mu(A\log k-1)+\log\log{k}\ge\log\frac{2s}\g+\log\log N=\log
4+\log\log N
\]
to guarantee (\ref{eq325}).
Again for large $k>k_8(\mu)$, we have $\mu+\log4\le\log\log k$,
and we only need to have
%
\begin{equation}\label{eq326}
k^{\mu A}\ge\log N.
\end{equation}
In summary, there exists a constant $k_9=k_9(q,s,\mu,A)$ such that
(\ref{eq324}) is valid
with $\tau_h=T_k(\delta)$ if $k>k_9$, $h=Ak\log k$ and $k$ satisfies
(\ref{eq315}), (\ref{eq30}) and
(\ref{eq326}) with $\mu\in(0,(q-1)/{2s})$.

\textit{Final step.} From (\ref{eq324}) and (\ref{eq316}) we learn
%
\begin{eqnarray}\label{eq327}
\bE_N T_k(\delta)&\le&
4(2-q)^{-1}k^{-s+2-q}+8k^{s(k-1)+3-q}N^{-1}\nonumber\\
&&{}+2\bigl({p(1-\delta)(1-\eta)}\bigr)^{-1}h^{1-\eta}\\
&&{}+2(1-\delta)^{-1}k^{-q/2}N^{{-q}/({2s(k-1)})}
h^3({ph})^{{q}/({p-1})},\nonumber
\end{eqnarray}
because $\tau_h=T_k(\delta)$. Condition (\ref{eq326}) combined with
(\ref{eq315}) and (\ref{eq30}) yield
%
\begin{equation}\label{eq328}
k^{\mu A}\ge\log N\ge(sk-s+2)\log k,\qquad 2k^{-s}\le1.
\end{equation}
For this to be plausible for large $k$, it suffices to have $\nu:=\mu
A>1$. Since $\mu\in(0,(q-1)/(2s))$, we pick some
\[
A>\frac{2s}{q-1}
\]
and select $\mu\in(A^{-1},(q-1)/(2s))$.
Since $h=Ak\log k$ and $p=2s(k-1)+1$, bound (\ref{eq327}) implies
\begin{eqnarray*}
\bE_N T_k(\delta)&\le& 4(2-q)^{-1}k^{-s+2-q}+8k^{s(k-1)+3-q}N^{-1}\\
&&{}+c_4(1-\delta)^{-1}k^{-\eta}(\log k)^{1-\eta}\\
&&{}+c_4(1-\delta
)^{-1}k^{3-q/2}(\log k)^3 N^{{-q}/({2s(k-1)})},
\end{eqnarray*}
because $({ph})^{{q}/({p-1})}$ is uniformly bounded in $k$. This
completes the proof of (\ref{eq31})
because (\ref{eq328}) is exactly (\ref{eq32}).
\end{pf*}

\section{Complete gelation}
\label{sec4}

This section is devoted to the proof of Theorem~\ref{th17}.
Lemma~\ref{lem4.1} below and Theorem~\ref{th16} are the main ingredients for the proof
of Theorem~\ref{th17}.

\begin{lemma}
\label{lem4.1}
Assume that $\alpha(m,n)\ge m^qn+n^qm$ for some $q>1$. Set
\[
\sigma=\inf\bigl\{t: L_{N/2}(t)>0\bigr\},\qquad
\hat \sigma=\min \bigl
\{\sigma,T_k(\delta) \bigr\}.
\]
Then
%
\begin{equation}
\label{eq4.1}
\bE_N (\sigma-\hat\sigma)\le4\delta^{-1}k^{1-q}.
\end{equation}
\end{lemma}

\begin{pf}
Define $K(\bL)=\sum_nL_n$.
By strong Markov property,
%
\begin{equation}
\label{eq4.2}
\bE_N K\bigl(\bL(\sigma)\bigr)=\bE_N K
\bigl(\bL(\hat\sigma)\bigr)+\bE_N \int^{\s
}_{\hat\s}
\cA K\bigl( \bL(t)\bigr)\,dt.
\end{equation}
If $\s>\hat\s$ and $\bL=\bL(t)$ for some $t\in({\hat\s},{\s
})$, then
\begin{eqnarray*}
-\cA K(\bL)&=& \frac 1{2N}\sum_{m,n}
\alpha(m,n)L_m \bigl(L_n-{\mathbh 1}(m=n) \bigr)
\\
&\ge&\frac1{2N}\sum_{m,n}n^qmL_m
\bigl(L_n-{\mathbh 1}(m=n) \bigr){\mathbh 1} (n\ge k)
\\
&=& \frac1{2N}\sum_{m,n}n^qmL_mL_n
{\mathbh 1}(n\ge k)-\frac1{2N}\sum_n
n^{q+1} L_n{\mathbh 1}(N/2\ge n\ge k)
\\
&\ge&\frac{1}{2N} \biggl(\sum_{n\ge
k}n^qL_n
\biggr) \biggl(\sum_{m}mL_m \biggr)-
\frac1{4}\sum_nn^qL_n{
\mathbh 1}(N/2\ge n\ge k)
\\
&=& \frac12\sum_{n\ge
k}n^qL_n-
\frac14\sum_{n\ge
k}n^qL_n=
\frac14 \sum_{n\ge k}n^qL_n \ge
\frac 14 N \delta k^{q-1}.
\end{eqnarray*}
From this and (\ref{eq4.2}) we deduce
\[
N\ge\bE_N K\bigl(\bL(\hat\sigma)\bigr) \ge\tfrac14 N \delta
k^{q-1}\bE_N (\sigma-\hat\sigma),
\]
as desired.
\end{pf}

\begin{pf*}{Proof of Theorem~\ref{th17}}
Let $\hat T=\hat T_A(1/2)$ be as in Theorem~\ref{th16}, with $A$ a positive
constant satisfying $A<q(2-q)^{-
1}(6-q)^{-1}$. Pick $\theta\in(0,\bar\eta)$ so that by (\ref{eq116}),
%
\begin{equation}
\label{eq4.3} \bE_N \hat T\le c_1 (\log
N)^{-\theta}
\end{equation}
for a constant $c_1$. Use Lemma~\ref{lem4.1} for $k=A\log N/\log\log N$ to
assert
\[
\bE_N \bigl(\sigma- \min\{\hat T,\sigma\} \bigr)\le
c_2 (\log N/\log\log N )^{1-q}
\]
for a
constant $c_2$. From this and (\ref{eq4.3}) we deduce
%
\begin{equation}
\label{eq4.4} \bE_N \sigma\le c_3 (\log N/\log\log N
)^{1-q}\vadjust{\goodbreak}
\end{equation}
for a constant $c_3$. Recall that at time $\s$,
we already have a particle of size at least $h=N/2$.
We mark one such particle and keep track of its size $\bar N(t)$ at
later times $t\ge\s$.
We also define an auxiliary process $(K(t)\dvtx t\ge\s_\ell)$ by the following
rules: $K(\s)=0$ and
$K$ increases by $1$ each time the marked particle coagulates with
another particle.
We would like to use this marked particle to produce a complete gelation.
Define the stopping time
\[
S_r=\inf \bigl\{t\dvtx\bar N(t)\ge r \bigr\}.
\]
Our goal is bounding $S_{r+1}-S_r$.
Note that if $S_{r+1}-S_r\neq0$, then $\bar N(t)=r$ for every $t\in
(S_r,S_{r+1})$, and
\[
K(S_{r+1})-K(S_r)=1,
\]
because any coagulation of the marked particle
results in $\bar N\ge r+1$.
As before we write $\cA$ for the generator of the augmented process
$\hat\bL(t)=(\bL(t),K(t))$ and
abuse the notation to write $K$ for the function that maps $\hat\bL$
to its second component $K$.
Note that if $\hat\bL=\hat\bL(t)$ for some $t\in(S_r,S_{r+1})$, then
\begin{eqnarray*}
\cA K(\hat\bL)&=& \frac1N\sum_m \alpha(r,m)\bigl[
L_m-{\mathbh 1}(m= r)\bigr]\ge\frac1N\sum
_m \bigl(r^q m+ m^qr\bigr)\bigl[
L_m-{\mathbh 1}(m= r)\bigr]
\\
&\ge&\frac1N\sum_m mr^q\bigl[
L_m-{\mathbh 1}(m= r)\bigr]= \biggl(1-\frac{r}N \biggr)
r^q.
\end{eqnarray*}
From this and strong Markov property
\[
1\ge\bE_N \bigl(K(S_{r+1})-K(S_r) \bigr)=
\bE_N \int^{S_{r+1}}_{S_r}\cA K\bigl(\hat
\bL(t)\bigr)\,dt,
\]
we deduce
\[
\bE_N (S_{r+1}-S_r)\le \biggl(1-
\frac{r}N \biggr)^{-1} r^{-q}.
\]
Summing this over $r$ yields
\begin{eqnarray*}
\bE_N (S_{N}-S_{h})&\le&\sum
_{r=h}^{N-1} \biggl(1-\frac{r}N
\biggr)^{-1} r^{-q} \le\frac N{h^q}\sum
_{r=h}^{N-1}(N-r)^{-1}
\\
&\le&\frac N{h^q} \bigl[\log(N-h)+1 \bigr]\le Nh^{-q}(1+
\log N).
\end{eqnarray*}
From this and (\ref{eq4.4}) we learn that if $\tilde\tau$ denotes
the time of the complete gelation,
then
\begin{eqnarray*}
\bE_N \tilde\tau&\le& c_3 (\log N/\log\log N
)^{1-q}+ 2^{-q}(1+\log N)N^{1-q}\\
&\le& c_4 (
\log N/\log\log N )^{1-q}.
\end{eqnarray*}
This completes the proof of (\ref{eq119}).
\end{pf*}

\section*{Acknowledgments}

The author thanks two anonymous referees for the references
\cite{FoL,FoG}
and \cite{N}. The author is also very grateful to one of the
referees who
spotted some errors on the earlier version of this paper; his comments
and suggestions have contributed greatly to the quality of this
article. Part of this work was done when the author was visiting Centre
Emile Borel of IHES and Universite de Paris 12. The author wishes to
thank IHES and Universite de Paris 12 for their invitation and
financial support.


%

\printaddresses


\begin{thebibliography}{16}

\bibitem{A2}
%
\begin{barticle}[mr]
\bauthor{\bsnm{Aldous},~\bfnm{David}\binits{D.}}
(\byear{1998}).
\btitle{Emergence of the giant component in special {M}arcus--{L}ushnikov
processes}.
\bjournal{Random Structures Algorithms}
\bvolume{12}
\bpages{179--196}.
\bid{doi={10.1002/(SICI)1098-2418(199803)12:2\&lt;179::AID-RSA2\&gt;3.0.CO;2-U},
issn={1042-9832}, mr={1637407}}
\bptok{imsref}%
\end{barticle}
%
\endbibitem


\bibitem{BC}
%
\begin{barticle}[mr]
\bauthor{\bsnm{Ball},~\bfnm{J.~M.}\binits{J.~M.}} \AND
\bauthor{\bsnm{Carr},~\bfnm{J.}\binits{J.}}
(\byear{1990}).
\btitle{The discrete coagulation--fragmentation equations: Existence,
uniqueness, and density conservation}.
\bjournal{J. Stat. Phys.}
\bvolume{61}
\bpages{203--234}.
\bid{issn={0022-4715}, mr={1084278}}
\bptok{imsref}%
\end{barticle}
%
\endbibitem

\bibitem{CdC}
%
\begin{barticle}[mr]
\bauthor{\bsnm{Carr},~\bfnm{J.}\binits{J.}} \AND\bauthor
{\bparticle{da}
\bsnm{Costa},~\bfnm{F.~P.}\binits{F.~P.}}
(\byear{1992}).
\btitle{Instantaneous gelation in coagulation dynamics}.
\bjournal{Z.~Angew. Math. Phys.}
\bvolume{43}
\bpages{974--983}.
\bid{doi={10.1007/BF00916423}, issn={0044-2275}, mr={1198671}}
\bptok{imsref}%
\end{barticle}
%
\endbibitem

\bibitem{EMP}
%
\begin{barticle}[mr]
\bauthor{\bsnm{Escobedo},~\bfnm{M.}\binits{M.}},
\bauthor{\bsnm{Mischler},~\bfnm{S.}\binits{S.}} \AND
\bauthor{\bsnm{Perthame},~\bfnm{B.}\binits{B.}}
(\byear{2002}).
\btitle{Gelation in coagulation and fragmentation models}.
\bjournal{Comm. Math. Phys.}
\bvolume{231}
\bpages{157--188}.
\bid{doi={10.1007/s00220-002-0680-9}, issn={0010-3616}, mr={1947695}}
\bptok{imsref}%
\end{barticle}
%
\endbibitem

\bibitem{Fl}
%
\begin{bbook}[auto:STB|2011/10/17|13:52:43]
\bauthor{\bsnm{Flory},~\bfnm{P.~J.}\binits{P.~J.}}
(\byear{1953}).
\btitle{Principle of Polymer Chemistry}.
\bpublisher{Cornell Univ. Press}, \baddress{Ithaca, NY}.
\bptok{imsref}%
\end{bbook}
%
\endbibitem

\bibitem{FoG}
%
\begin{barticle}[mr]
\bauthor{\bsnm{Fournier},~\bfnm{Nicolas}\binits{N.}} \AND
\bauthor{\bsnm{Giet},~\bfnm{Jean-S{\'e}bastien}\binits{J.-S.}}
(\byear{2004}).
\btitle{Convergence of the {M}arcus--{L}ushnikov process}.
\bjournal{Methodol. Comput. Appl. Probab.}
\bvolume{6}
\bpages{219--231}.
\bid{doi={10.1023/B:MCAP.0000017714.56667.67}, issn={1387-5841}, mr={2035293}}
\bptok{imsref}%
\end{barticle}
%
\endbibitem

\bibitem{FoL}
%
\begin{barticle}[mr]
\bauthor{\bsnm{Fournier},~\bfnm{Nicolas}\binits{N.}} \AND
\bauthor{\bsnm{Lauren{\c{c}}ot},~\bfnm{Philippe}\binits{P.}}
(\byear{2009}).
\btitle{Marcus--{L}ushnikov processes, {S}moluchowski's and {F}lory's models}.
\bjournal{Stochastic Process. Appl.}
\bvolume{119}
\bpages{167--189}.
\bid{doi={10.1016/j.spa.2008.02.003}, issn={0304-4149}, mr={2485023}}
\bptok{imsref}%
\end{barticle}
%
\endbibitem


\bibitem{HR1}
%
\begin{barticle}[mr]
\bauthor{\bsnm{Hammond},~\bfnm{Alan}\binits{A.}} \AND
\bauthor{\bsnm{Rezakhanlou},~\bfnm{Fraydoun}\binits{F.}}
(\byear{2007}).
\btitle{Moment bounds for the {S}moluchowski equation and their consequences}.
\bjournal{Comm. Math. Phys.}
\bvolume{276}
\bpages{645--670}.
\bid{doi={10.1007/s00220-007-0304-5}, issn={0010-3616}, mr={2350433}}
\bptnote{check year}%
\bptok{imsref}%
\end{barticle}
%
\endbibitem

\bibitem{HR2}
%
\begin{barticle}[mr]
\bauthor{\bsnm{Hammond},~\bfnm{Alan}\binits{A.}} \AND
\bauthor{\bsnm{Rezakhanlou},~\bfnm{Fraydoun}\binits{F.}}
(\byear{2007}).
\btitle{The kinetic limit of a system of coagulating {B}rownian particles}.
\bjournal{Arch. Ration. Mech. Anal.}
\bvolume{185}
\bpages{1--67}.
\bid{doi={10.1007/s00205-006-0033-5}, issn={0003-9527}, mr={2308858}}
\bptok{imsref}%
\end{barticle}
%
\endbibitem

\bibitem{J1}
%
\begin{barticle}[mr]
\bauthor{\bsnm{Jeon},~\bfnm{Intae}\binits{I.}}
(\byear{1998}).
\btitle{Existence of gelling solutions for coagulation--fragmentation
equations}.
\bjournal{Comm. Math. Phys.}
\bvolume{194}
\bpages{541--567}.
\bid{doi={10.1007/s002200050368}, issn={0010-3616}, mr={1631473}}
\bptok{imsref}%
\end{barticle}
%
\endbibitem

\bibitem{J2}
%
\begin{barticle}[mr]
\bauthor{\bsnm{Jeon},~\bfnm{Intae}\binits{I.}}
(\byear{1999}).
\btitle{Spouge's conjecture on complete and instantaneous gelation}.
\bjournal{J. Stat. Phys.}
\bvolume{96}
\bpages{1049--1070}.
\bid{doi={10.1023/A:1004640317274}, issn={0022-4715}, mr={1722986}}
\bptok{imsref}%
\end{barticle}
%
\endbibitem

\bibitem{N}
%
\begin{barticle}[mr]
\bauthor{\bsnm{Norris},~\bfnm{James~R.}\binits{J.~R.}}
(\byear{1999}).
\btitle{Smoluchowski's coagulation equation: Uniqueness, nonuniqueness
and a
hydrodynamic limit for the stochastic coalescent}.
\bjournal{Ann. Appl. Probab.}
\bvolume{9}
\bpages{78--109}.
\bid{doi={10.1214/aoap/1029962598}, issn={1050-5164}, mr={1682596}}
\bptok{imsref}%
\end{barticle}
%
\endbibitem

\bibitem{R}
%
\begin{barticle}[mr]
\bauthor{\bsnm{Rezakhanlou},~\bfnm{F.}\binits{F.}}
(\byear{2006}).
\btitle{The coagulating {B}rownian particles and {S}moluchowski's equation}.
\bjournal{Markov Process. Related Fields}
\bvolume{12}
\bpages{425--445}.
\bid{issn={1024-2953}, mr={2249642}}
\bptok{imsref}%
\end{barticle}
%
\endbibitem

\end{thebibliography}
\end{document}